\documentclass[11pt,a4paper]{article}

\usepackage[utf8]{inputenc}

\usepackage{fullpage}
\usepackage{hyperref}

\usepackage{graphicx}

\usepackage{amsmath}
\usepackage{amsfonts}
\usepackage{amsthm}

\usepackage{enumerate}

\newcommand{\ie}{i.e.}
\newcommand{\eg}{e.g.}

\theoremstyle{definition} 
\newtheorem{definition}{Definition}[section] 
\theoremstyle{remark} 
\newtheorem{remark}{Remark} 

\newcommand{\eqdef}{\mathrel{:=}}
\newcommand{\reveqdef}{\mathrel{=:}}

\newcommand{\tr}[1]{#1^\mathit{T}}
\newcommand{\conj}[1]{\overline{#1}}

\newcommand{\Rbb}{\mathbb{R}}
\newcommand{\Sbb}{\mathbb{S}}

\newcommand{\tgspace}[2]{T_{#1}{#2}}
\newcommand{\dirder}[3]{D#1(#2)[#3]}
\newcommand{\Euclgrad}[3]{\nabla_{#3}#1(#2)}

\newcommand{\myfigurename}{Fig.}

\newcommand{\params}{\lambda}
\newcommand{\Params}{\Lambda}
\newcommand{\PRC}{q_\epsilon}
\newcommand{\iPRCx}{q_x}
\newcommand{\iPRCu}{q}
\newcommand{\chifield}{V}
\newcommand{\PLC}{\psi}

\begin{document}

\title{Sensitivity analysis of circadian entrainment \\ in the space of phase response curves%
	\thanks{This chapter presents research results of the Belgian Network DYSCO (Dynamical Systems, Control, and Optimization), funded by the Interuniversity Attraction Poles Programme, initiated by the Belgian State, Science Policy Office. The scientific responsibility rests with its authors. P.~Sacr\'{e} is a Research Fellow with the Belgian Fund for Scientific Research (F.R.S.-FNRS).}
	}

\author{Pierre~Sacr\'{e} and~Rodolphe~Sepulchre
	\thanks{P.~Sacr\'{e} and R.~Sepulchre are with the Department of Electrical Engineering and Computer Science (Montefiore Institute), University of Li\`{e}ge, 4000 Li\`{e}ge, Belgium (e-mails: pierre.sacre@ulg.ac.be, r.sepulchre@ulg.ac.be).}
	}
	
\date{}
	
\maketitle

\begin{abstract}
	Sensitivity analysis is a classical and fundamental tool to evaluate the role of a given parameter in a given system characteristic. Because the phase response curve is a fundamental input--output characteristic of oscillators, we developed a sensitivity analysis for oscillator models in the space of phase response curves. The proposed tool can be applied to high-dimensional  oscillator models without facing the curse of dimensionality obstacle associated with numerical exploration of the parameter space. Application of this tool to a state-of-the-art model of circadian rhythms suggests that it can be useful and instrumental to biological investigations.
\end{abstract}

\section{Introduction} 
\label{sec:introduction}

\begin{figure}
	\centering
	\includegraphics{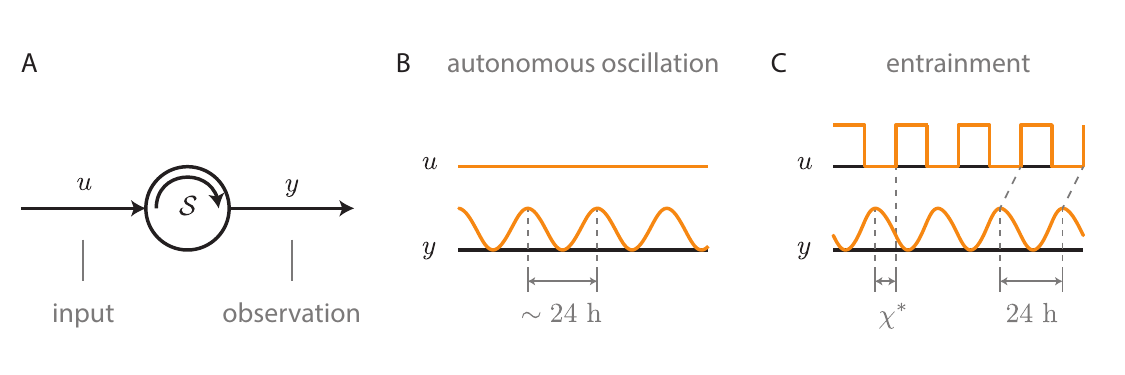}
	\caption{(A)~Circadian oscillators are viewed as open dynamical systems with input~$u$ and output~$y$. (B)~The unforced system exhibits autonomous rhythms that occur with a period close to $24$~hours. (C)~The periodically forced system adapts the organism rhythms through entrainment ($1\mathord{:}1$ phase-locking) with the 24-hours stimulus associated with earth's rotation.}
	\label{fig:circ-osc}
\end{figure}

Circadian entrainment is a biological process at the core of most living organisms which need to adapt their physiological activity to the 24 hours environmental cycle associated with earth's rotation (\eg{} variations in light or temperature condition). 
This process relies on the robust interaction between an autonomous molecular oscillator and its environment (\myfigurename~\ref{fig:circ-osc}A). Experimental observations have shown that the system is capable to exhibit oscillations with a period close to 24 hours in constant environmental condition (unforced system, \myfigurename~\ref{fig:circ-osc}B) and to lock its oscillations (in frequency and phase) to an environmental cue with a period equal to 24 hours (periodically forced system, \myfigurename~\ref{fig:circ-osc}C). This locking phenomenon is often called (circadian) entrainment \cite{Pittendrigh:1981wa}.
Moreover, this biological process is known to be very robust, that is, it maintains its performance (its period and its locking) despite internal or external perturbations (\eg{}~genetic mutations, molecular noise, variability of the environmental condition,~etc.).

With recent experimental advances in  biology, the molecular bases of circadian rhythms has been increasingly unfolded in various organisms. In most eukaryotic organisms (\eg{}~fungus, fly, or mouse), the core mechanism relies on analogous interacting positive and negative feedback loops with several minor alterations~\cite{Hastings:2000ht}. However, even though the architecture of those biological clocks is better known, the specific design and robustness mechanisms implemented in those architectures remain unknown \cite{Stelling:2004do,Novak:2008eg}.

Starting with the pioneering work of Winfree~\cite{Winfree:1967vf,Winfree:1980ue}, the Phase Response Curve (PRC)  has emerged as a fundamental input--output characteristic of oscillators. 
Analogously to the static gain of a transfer function, the PRC measures a steady-state (asymptotic) property of the system response induced by an impulsive input. For the static gain, the measured property is the integral of the response; for the PRC, the measured property is the phase shift between the unperturbed and perturbed responses. 
Because of the periodic nature of the steady-state, this phase shift depends on the phase at which the system receives the impulsive input. The PRC is thus a curve rather than a scalar.
In many situations, the PRC can be determined experimentally and provides unique data for the model identification of the oscillator. Likewise, numerical methods exist to compute the PRC from a state-space model of the oscillator. Finally, the PRC contains the fundamental mathematical information required to reduce a $n$-dimensional state-space model to the one-dimensional (phase) center manifold of a hyperbolic limit cycle.

In this chapter, we review (local) sensitivity tools that provide numerical and mathematical grounds to the robustness analysis of oscillator state-space models in  connection with experimentally available observations like the PRC or the period. We then illustrate how these tools can be used to make physiologically relevant predictions from  mathematical models of  circadian rhythms. We apply our sensitivity analysis to a state-of-the-art model \cite{Leloup:2003cp} of 16~states and 52~parameters and exploit the results to extract the parameters and circuits that determine the robustness of entrainment.

The \emph{local} proposed approach is systematic and computationally tractable. It provides a rapid screening of all parameters, even in high-dimensional models with a large number of parameters. It complements \emph{nonlocal} analyses often used to assess the robustness of parameters, such as bifurcation analysis \cite{Leloup:2004co} or parameter space exploration \cite{Stelling:2004do,Hafner:2009eb}.

The chapter is organized as follows.
Section~\ref{sec:open_oscillator_models} reviews the notion of PRCs characterizing the input--output behavior of an oscillator model in the neighborhood of a stable limit cycle.
Section~\ref{sec:sensitivity} develops the sensitivity analysis for oscillators in terms of the sensitivity of its periodic orbit, its PRC, and its entrainment (phase-locking).
Section~\ref{sec:scalar_robustness_measures_for_oscillators} provides scalar robustness measures based on this sensitivity analysis.
Section~\ref{sec:illustations} illustrates how those tools permit to address system-theoretic questions meaningful for the robustness analysis of circadian entrainment.


\section{Open oscillator models: from state-space to phase models} 
\label{sec:open_oscillator_models}

In this section, we provide a short introduction to oscillators viewed as open dynamical systems, that is, as dynamical systems that interact with their environment \cite{Sepulchre:2006vk}.
We first recall basic definitions about stable periodic orbits in $n$-dimensional state-space models (see~\cite{Farkas:1994uq,Khalil:2002wj} for details).
We then introduce (finite and infinitesimal) phase response curves as fundamental input--output mathematical information required for the model reduction.
We finally summarize the standard phase reduction procedure which concentrates the phase behavior information of $n$-dimensional state-space models into one-dimensional phase models characterized by its angular frequency, its PRC, and a measurement map (see~\cite{Kuramoto:1984wo,Hoppensteadt:1997tp} for details).

\subsection{State-space models: periodic orbits and phase maps}

We consider open dynamical systems described by nonlinear (single-input and single-output\footnote{For presentation convenience, we consider single-input and single-output systems. All developments are easily generalizable to multiple-input and multiple-output systems.}) time-invariant state-space models
\begin{subequations} \label{eq:nlsys}
	\begin{align}
		\dot{x} & = f(x) + g(x) u,  & x & \in\Rbb^n, \; u\in\Rbb, \label{eq:nlsys_dx} \\
			 y  & = h(x), 			& y & \in\Rbb, \label{eq:nlsys_y}
	\end{align}	
\end{subequations}
where the \emph{vector fields}~$f$ and~$g$, and the \emph{measurement map}~$h$ support all the usual smoothness conditions that are necessary for existence and uniqueness of solutions. We denote by~$\phi(\cdot,x_0,u)$ the solution to the initial value problem~\eqref{eq:nlsys_dx} from the initial condition~$x_0\in \Rbb^n$ at time $0$, that is, $\phi(0,x_0,u) = x_0$. 

\begin{figure}
	\centering
	\includegraphics{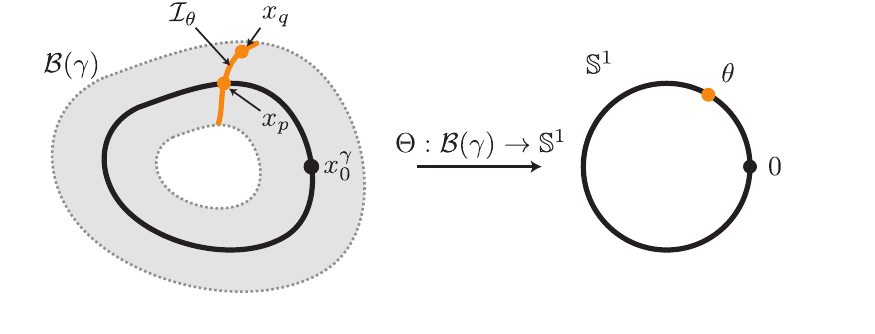}
	\caption{The asymptotic phase map~$\Theta:\mathcal{B}(\gamma)\rightarrow\Sbb^1$ associates with each point~$x_q$ in the basin~$\mathcal{B}(\gamma)$ a scalar phase~$\Theta(x_q)=\theta$ on the unit circle~$\Sbb^1$ such that~$\lim_{t\rightarrow+\infty} \left\|\phi(t,x_q,0) - \phi(t,x_p,0)\right\|_2 = 0$ with~$x_p = x^\gamma(\theta/\omega)$.}
	\label{fig:notations}
\end{figure}

An oscillator is an open dynamical system whose zero-input steady-state behavior is periodic rather than constant. Formally, we assume that the zero-input system~$\dot{x}=f(x)$ admits a \emph{locally hyperbolic stable periodic orbit}~$\gamma\subseteq\Rbb^n$ with period~$T$ (and corresponding angular frequency~$\omega=2\pi/T$). 
Picking an initial condition~$x^\gamma_0$ on the periodic orbit~$\gamma$, this latter is described by the (nonconstant) $T$-periodic trajectory~$\phi(\cdot,x^\gamma_0,0)\reveqdef x^\gamma(\cdot)$, such that $x^\gamma(\cdot)=x^\gamma(\cdot+T)$.
The \emph{basin of attraction} of~$\gamma$ is the maximal open set from which the periodic orbit~$\gamma$ attracts. (Main notations are illustrated on \myfigurename~\ref{fig:notations}.)

Since the periodic orbit~$\gamma$ is a one-dimensional manifold in~$\Rbb^n$, it is homeomorphic to the unit circle~$\mathbb{S}^1$. It is thus naturally parametrized in terms of a single scalar \emph{phase}.
The smooth bijective \emph{phase map} $\varTheta:\gamma\rightarrow\Sbb^1$ associates with each point $x_p$ on the periodic orbit $\gamma$ its phase $\varTheta(x_p)\reveqdef\vartheta_p$ on the unit circle $\Sbb^1$, such that,
\begin{equation}
	x_p - x^\gamma(\vartheta_p/\omega) = 0.
\end{equation}
This mapping is constructed such that the image of the reference point $x_0^\gamma$ is equal to~$0$ (\ie{} $\varTheta(x^\gamma_0)=0$) and the progression along the periodic orbit (in absence of perturbation) produces a constant increase in $\vartheta$.
The \emph{phase variable} $\vartheta:\Rbb_{\geq 0} \rightarrow \Sbb^1$ is defined along each zero-input trajectory $\phi(\cdot,x_0,0)$ starting from a point $x_0$ on the periodic orbit $\gamma$, as $\vartheta(t) \eqdef \varTheta(\phi(t,x_0,0))$ for all times $t \geq 0$. The phase dynamics are thus given by $\dot{\vartheta} = \omega$.

For hyperbolic stable periodic orbit, the notion of phase can be extended to any point~$x_q$ in the basin~$\mathcal{B}(\gamma)$ by defining the concept of \emph{asymptotic phase}. 
The \emph{asymptotic phase map}~$\Theta : \mathcal{B}(\gamma) \rightarrow \mathbb{S}^1$ associates with each point~$x_q$ in the basin~$\mathcal{B}(\gamma)$ its asymptotic phase~$\Theta(x_q) \reveqdef \theta_q$ on the unit circle $\Sbb^1$, such that, 
\begin{equation}
	\lim_{t\rightarrow+\infty} \left\|\phi(t,x_q,0) - \phi(t,x^\gamma(\theta_q/\omega),0)\right\|_2 = 0.
\end{equation}
Again, this mapping is constructed such that the image of~$x^\gamma_0$ is equal to~$0$ and such that the progression along any orbit in~$\mathcal{B}(\gamma)$ (in absence of perturbation) produces a constant increase in~$\theta$.
The \emph{asymptotic phase variable} $\theta:\Rbb_{\geq 0} \rightarrow \Sbb^1$ is defined along each zero-input trajectory $\phi(\cdot,x_0,0)$ starting from a point $x_0$ in the basin of attraction of~$\gamma$, as $\theta(t) \eqdef \Theta(\phi(t,x_0,0))$ for all times $t \geq 0$. The asymptotic phase dynamics are thus given by $\dot{\theta} = \omega$.
 
The notion of asymptotic phase variable can be extended to a nonzero-input trajectory $\phi(\cdot,x_0,u)$ provided that its stays in the basin of attraction of $\gamma$. In this case, the asymptotic phase variable is defined as $ \theta(t) \eqdef \Theta(\phi(t,x_0,u))$ for all times $t \ge 0$. 
Thus the variable $\theta(t_*)$ at an instant $t_* \ge 0$ evaluates the asymptotic phase of the point $\phi(t_*,x_0 ,u)$ such that
\begin{equation}
	\lim_{t\to+\infty} \left\|\phi(t,\phi(t_*,x_0,u),0)-\phi(t,x^\gamma(\theta(t_*)/\omega),0)\right\| = 0.
\end{equation}
The asymptotic phase dynamics in the case of a nonzero input is often hard to derive. 

For presentation convenience, we introduce the map $\tilde{x}^\gamma:\Sbb^1\rightarrow\gamma$ which associates to each phase $\theta$ a point $\phi(\theta/\omega,x_0^\gamma,0)=\tilde{x}^\gamma(\theta)$ on the periodic orbit. This map corresponds to a reparametrization of the periodic solution~$x^\gamma(\cdot)$.

The $2\pi$-periodic steady-state solution $\tilde{x}^\gamma(\cdot)$ and the angular frequency $\omega$ can be calculated by solving the boundary value problem \cite{Ascher:1988ty,Seydel:2010tw}
\begin{subequations} \label{eq:bvp_x}
	\begin{align}
		(\tilde{x}^\gamma)'(\theta) - \frac{1}{\omega} f(\tilde{x}^\gamma(\theta)) & = 0 \label{eq:bvp_x_dx} \\
		\tilde{x}^\gamma(2\pi) - \tilde{x}^\gamma(0) & = 0 \label{eq:bvp_x_per}\\
		\PLC(\tilde{x}^\gamma(0),\omega)  & = 0 \label{eq:bvp_x_phase}
	\end{align}	
\end{subequations}
(where the prime $\cdot'$ denotes the derivative with respect to~$\theta$). The boundary conditions are given by the periodicity condition~\eqref{eq:bvp_x_per} which ensures the periodicity of the map~$\tilde{x}^\gamma(\cdot)$ and the phase condition~\eqref{eq:bvp_x_phase} which anchors a reference position~$\tilde{x}^\gamma(0)$ along the periodic orbit. The phase condition $\PLC:\Rbb^n \times \Rbb_{>0} \rightarrow \Rbb$ is chosen such that it defines an isolated point on the periodic orbit (see~\cite{Seydel:2010tw} for details).
Numerical algorithms to solve this boundary value problem are reviewed in \cite[Appendix]{Sacre:2012uz}.
 
\subsection{Phase response curves: local information about the phase map}

For many oscillators, the structure of the asymptotic phase map is very complex. This often makes its analytical computation impossible and even its numerical computation intractable (or at least very expensive, in particular for high-dimensional oscillator models). However, in many situations, the global knowledge of the asymptotic phase map is not required to study oscillator dynamics. Instead, it is sufficient to consider a local phase information also known as the phase response curve.

Starting with the pioneering work of Winfree \cite{Winfree:1967vf,Winfree:1980ue}, the phase response curve of an oscillator has proven a useful input--output tool to study oscillator dynamics. It indicates how the timing of inputs affects the timing (steady-state phase shift) of oscillators. Phase response curves are directly related to asymptotic phase maps but capture only partial (local) information about them.

\begin{definition}
	The \emph{Phase Response Curve (PRC)} corresponding to an impulsive input of finite amplitude $\epsilon$ (\ie{} $u(\cdot) \eqdef \epsilon \delta(\cdot)$ where $\delta(\cdot)$ is the Dirac delta function) is the map $\PRC:\mathbb{S}^1\rightarrow(-\pi,\pi]$ defined as
	\begin{equation}
		\PRC(\theta) \eqdef \Delta \Theta(\tilde{x}^\gamma(\theta)) = \lim_{t\rightarrow 0^+}\underbrace{\Theta( \phi(t,\tilde{x}^\gamma(\theta),\epsilon \delta(\cdot)) )}_{\text{post-stimulus phase}}  -  \underbrace{\Theta( \phi(t,\tilde{x}^\gamma(\theta),0) )}_{\text{pre-stimulus phase}}.
	\end{equation}
	It associates with each point on the periodic orbit (parametrized by its phase $\theta$) the phase shift induced by the input.
\end{definition}

In many situations, the PRC can be determined experientially (in particular for circadian rhythms). Moreover, it can be computed numerically by simulating the nonlinear state-space model and comparing the asymptotic phase shift between perturbed and unperturbed trajectories.

A mathematically more abstract---yet very useful---tool is the infinitesimal phase response curve. It records essentially the same information as the finite phase response curve but for infinitesimally small Dirac delta input ($\epsilon \ll 1$). 
\begin{definition}
	The (input) \emph{infinitesimal Phase Response Curve (iPRC)} is the map $\iPRCu:\mathbb{S}^1\rightarrow\mathbb{R}$ defined as the directional derivative 
	\begin{equation} \label{eq:q_u}
		\iPRCu(\theta) \eqdef \dirder{\Theta}{\tilde{x}^\gamma(\theta)}{g(\tilde{x}^\gamma(\theta))}
	\end{equation}
	where
	\begin{equation}
		\dirder{\Theta}{x}{\vec{\eta}} \eqdef \lim_{\epsilon\rightarrow0}\frac{\Theta(x+\epsilon\vec{\eta}) - \Theta(x)}{\epsilon}.
	\end{equation}
	The directional derivative can be computed as the inner product
	\begin{equation}
		\dirder{\Theta}{x}{g(x)} = \langle \Euclgrad{\Theta}{x}{x} , g(x) \rangle
	\end{equation}
	where $\Euclgrad{\Theta}{x}{x}$ is the gradient of $\Theta$ at $x$. The map $\iPRCx:\Sbb^1\rightarrow\Rbb^n:\theta\mapsto\Euclgrad{\Theta}{\tilde{x}^\gamma(\theta)}{x}\reveqdef \iPRCx(\theta)$ is known as the state infinitesimal phase response curve.
\end{definition}

The (state) iPRC $\iPRCx(\cdot)$ can be calculated  by solving the boundary value problem \cite{Kuramoto:1984wo,Malkin:1949to,Malkin:1956wq,Neu:1979cj,Ermentrout:1984jb}
\begin{subequations} \label{eq:bvp_q}
	\begin{align}
		\iPRCx'(\theta) + \frac{1}{\omega} \tr{f_x(\tilde{x}^\gamma(\theta))} \iPRCx(\theta) & = 0 \label{eq:bvp_q_a} \\
		\iPRCx(2\pi) - \iPRCx(0) & = 0 \label{eq:bvp_q_b} \\
		\langle \iPRCx(\theta) , f(\tilde{x}^\gamma(\theta)) \rangle - \omega & = 0 \label{eq:bvp_q_c}
	\end{align}	
\end{subequations}
(where the notation $\tr{A}$ stands for the transpose of the matrix $A$). The boundary condition~\eqref{eq:bvp_q_b} imposes the periodicity of~$\iPRCx(\cdot)$ and the normalization condition~\eqref{eq:bvp_q_c}~ensures a linear increase at rate $\omega$ of the phase variable~$\theta$ along zero-input trajectories.
Numerical methods to solve this boundary value problem as a by-product of the periodic orbit computation are presented in \cite[Appendix]{Sacre:2012uz}.

\begin{remark}
	For small values of $\epsilon$ (\ie{} $\epsilon \ll 1$), the PRC for impulsive input of finite amplitude is well approximated by the iPRC, that is, $\PRC(\cdot) = \epsilon \iPRCu(\cdot) + \mathcal{O}(\epsilon^2)$.
\end{remark}

\subsection{Phase models: entrainment} \label{sec:phasemodel}

In the weak perturbation limit, that is, for small inputs
\begin{equation}
	u(t) = \epsilon \bar{u}(t), \quad \epsilon \ll 1, \quad |\bar{u}(t)| \leq 1 \quad \text{for all $t$},
\end{equation}
any solution $\phi(t,x_0,u)$ of the oscillator model which starts in the neighborhood of the hyperbolic stable periodic orbit~$\gamma$ stays in its neighborhood. The $n$-dimensional state-space model can thus be approximated by a one-dimensional (continuous-time) phase model \cite{Kuramoto:1984wo,Malkin:1949to,Malkin:1956wq,Neu:1979cj,Ermentrout:1984jb}
\begin{subequations}
	\begin{align} \label{eq:weak_coupling}
		\dot{\theta} & = \omega + \epsilon \iPRCu(\theta) \bar{u}(t) & \theta & \in\Sbb, \; \bar{u}\in\Rbb, \\
		y            & = \tilde{h}(\theta) & y & \in\Rbb.
	\end{align}	
\end{subequations}
The phase model is fully characterized by its angular frequency $\omega>0$, its infinitesimal phase response map $\iPRCu:\mathbb{S}^1 \rightarrow \mathbb{R}$, and its measurement map $\tilde{h}:\mathbb{S}^1 \rightarrow \mathbb{R}$.

\begin{figure}
	\centering
	\includegraphics{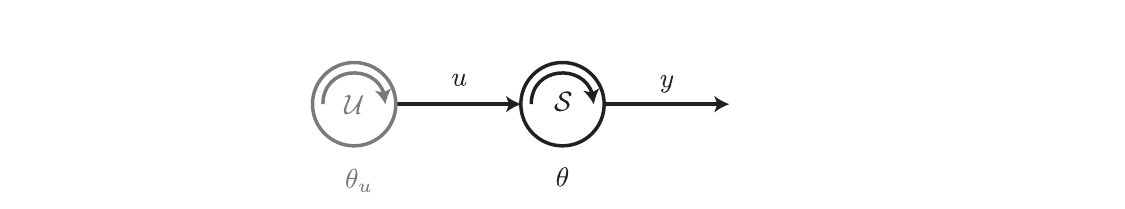}
	\caption{Entrainment is studied by applying weakly connected oscillator theory to the feedforward interconnection between an artificial oscillator generating the input and the actual oscillator.}
	\label{fig:entrainment}
\end{figure}

To study entrainment through weak coupling, we can apply weakly connected oscillator theory \cite[Chapter 9]{Hoppensteadt:1997tp} by considering the input $u(t)$ as generated by an artificial oscillator described by the trivial phase model 
$\dot{\theta}_u = \omega_u, \; y_u = \tilde{h}_u(\theta_u)$,
where we denote by $\omega_u$ the input angular frequency and we choose the artificial oscillator output map $\tilde{h}_u$ such that $y_u(t) = \bar{u}(t)$ for all times $t \geq 0$. Moreover, the network interconnection in this case is a feedforward interconnection from the artificial oscillator generating the input to the studied oscillator (see \myfigurename~\ref{fig:entrainment}).

The interconnected phase dynamics are thus given by
\begin{subequations} \label{eq:phase_dyn}
	\begin{align}
		\dot{\theta}_u & = \omega_u \\
		\dot{\theta}   & = \omega + \epsilon \iPRCu(\theta) \tilde{h}_u(\theta_u).
	\end{align}	
\end{subequations}
Following the weakly connected oscillator theory, we decompose the angular frequencies as $\omega = \Omega + \Delta$ and $\omega_u = \Omega_u + \Delta_u$ with $\Omega - \Omega_u = 0$, and the phase variables  as  $\theta=\Omega t+\varphi$ and $\theta_u=\Omega_u t+\varphi_u$ where $\varphi$ and $\varphi_u$ are slow phase deviations from the fast oscillations $\Omega t$ and $\Omega_u t$. The phase deviation dynamics are given by
\begin{subequations}
	\begin{align}
		\dot{\varphi}_u & = \Delta_u \\
		\dot{\varphi}   & = \Delta + \epsilon \iPRCu(\Omega t + \varphi) \tilde{h}_u(\Omega_u t + \varphi_u) .
	\end{align}	
\end{subequations}
Assuming that $\Delta,\Delta_u,\epsilon \ll 1$, standard averaging techniques yield
\begin{subequations} \label{eq:phase_dev_dyn}
	\begin{align}
		\dot{\varphi}_u & = \Delta_u \\
		\dot{\varphi}   & = \Delta + \epsilon \Gamma(\varphi - \varphi_u)
	\end{align}	
\end{subequations}
where the coupling function is given by
\begin{equation} \label{eq:coupl_fun}
	\Gamma(\varphi - \varphi_u) = \lim_{\tilde{T}\rightarrow+\infty}\frac{1}{\tilde{T}} \int_0^{\tilde{T}} \iPRCu(\Omega t + \varphi - \varphi_u) \tilde{h}_u(\Omega_u t ) dt .
\end{equation}
Introducing the phase difference $\chi = \varphi - \varphi_u$, we have
\begin{equation} \label{eq:phase_diff_dyn}
	\dot{\chi} =  \Delta - \Delta_u + \Gamma(\chi) \reveqdef \chifield(\chi).
\end{equation}
A stable equilibrium $\chi^*$ of \eqref{eq:phase_diff_dyn}, that is,
\begin{equation}
	\chi^* \in \Sbb^1 : \chifield(\chi^*) = 0 \quad \text{and} \quad \chifield'(\chi^*) < 0
\end{equation}
correspond to a stable  $1\mathord{:}1$~phase-locking behavior (or entrainment) for~\eqref{eq:phase_dyn}, that is,
\begin{equation}
	\theta(t) - \theta_u(t) = \chi^* \quad \text{for all times $t$}.
\end{equation}

\bigskip

In this section, we saw that a state-space oscillator model may be reduced to a phase model characterized by its angular frequency (or period) and its (infinitesimal) phase response curve. In addition, phase models are very useful to study entrainment. It is then very natural to study the sensitivity of oscillators with an emphasis on those characteristics.


\section{Sensitivity analysis for oscillators} 
\label{sec:sensitivity}

Sensitivity analysis for oscillators has been widely studied in terms of sensitivity analysis of periodic orbits~\cite{Kramer:1984jk,Rosenwasser:1999tw,Ingalls:2004wr,Wilkins:2009kq}. Because the phase response curve is an important oscillator characteristic, we recently proposed a sensitivity analysis of oscillator models in the space of phase response curves \cite{Sacre:2012uz}.  Moreover, the sensitivity analysis in the space of PRC can be exploited to predict the sensitivity of the entrainment.

We summarize those developments for nonlinear time-invariant state-space models with one parameter\footnote{For presentation convenience, we consider systems with a one-dimensional parameter space. All developments are easily generalizable to systems with a $q$-dimensional parameter space.}
\begin{subequations}
	\begin{align}
		\dot{x} & = f(x,\params) + g(x,\params) u \\
			 y  & = h(x,\params)
	\end{align}	
\end{subequations}
where the constant parameter $\params$ belongs to $\Rbb$.

\subsection{Sensitivity analysis of a periodic orbit}

The periodic orbit $\gamma$ of an oscillator model is characterized by its angular frequency~$\omega$ which measures the `speed' of a solution along the orbit and by the $2\pi$-periodic steady-state solution~$\tilde{x}^\gamma(\cdot)$ which describes the locus of this orbit in the state space. The sensitivity of both characteristics is important.

Given a nominal parameter value $\params_0$, the sensitivity of the angular frequency is the scalar $S_\omega \in \Rbb$ defined as
\begin{equation}
	S_\omega \eqdef \left.\frac{d\omega}{d\params}\right|_{\params_0} = \lim_{h \rightarrow 0} \frac{\left.\omega\right|_{\params_0 + h} - \left.\omega\right|_{\params_0}}{h}
\end{equation}
where the notation $\left.\star\right|_{\params}$ emphasizes the parameter value $\params$ at which the model characteristic~$\star$ is evaluated. 
Likewise, the sensitivity of the $2\pi$-periodic steady-state solution is the $2\pi$-periodic function $Z_{\tilde{x}}:\Sbb^1 \rightarrow \Rbb^{n}$ defined as
\begin{equation}
	Z_{\tilde{x}}(\cdot) \eqdef \left.\frac{d\tilde{x}^\gamma}{d\params}(\cdot)\right|_{\params_0} = \lim_{h \rightarrow 0} \frac{\left.\tilde{x}^\gamma(\cdot)\right|_{\params_0 + h} - \left.\tilde{x}^\gamma(\cdot)\right|_{\params_0}}{h}
\end{equation}
where the explicit dependence of the $2\pi$-periodic steady-state solution in $\params$ is given by
\begin{equation}
	\tilde{x}^\gamma(\cdot)|_{\params} = \left.\phi(\cdot/\omega|_{\params},x_0^\gamma|_{\params},0)\right|_{\params} . 	
\end{equation}
From~\eqref{eq:bvp_x}, we have, taking derivatives with respect to~$\params$,
\begin{subequations} \label{eq:bvp_Zx}
	\begin{align}
		Z_{\tilde{x}}'(\theta) - \frac{1}{\omega}A(\theta) Z_{\tilde{x}}(\theta) + \frac{1}{\omega^2}\tilde{v}(\theta)S_{\omega} - \frac{1}{\omega} b(\theta) & = 0 \\
		Z_{\tilde{x}}(2\pi) - Z_{\tilde{x}}(0) & = 0 \\
		\PLC_x Z_{\tilde{x}}(0) + \PLC_\omega S_{\omega} + \PLC_{\params} & = 0 
	\end{align}
\end{subequations}
where we use the following short notations
\begin{align}
	A(\cdot)        & \eqdef \frac{\partial f}{\partial x}(\tilde{x}^\gamma(\cdot),\params_0), & \PLC_x & \eqdef \frac{\partial \PLC}{\partial x}(x_0^\gamma,\omega,\params_0), \\
	b(\cdot)        & \eqdef \frac{\partial f}{\partial \params}(\tilde{x}^\gamma(\cdot),\params_0), & \PLC_\omega & \eqdef \frac{\partial \PLC}{\partial \omega}(x_0^\gamma,\omega,\params_0), \\
	\tilde{v}(\cdot)& \eqdef f(\tilde{x}^\gamma(\cdot),\params_0), & \PLC_\params & \eqdef \frac{\partial \PLC}{\partial \params}(x_0^\gamma,\omega,\params_0).
\end{align}

\begin{remark} \label{rem:period_sens}
	In the literature, the sensitivity of the period is often used instead of the sensitivity of the angular frequency. It is the scalar $S_T\in\Rbb$ defined as 
	\begin{equation}
		S_T \eqdef \left.\frac{dT}{d\params}\right|_{\params_0} = \lim_{h \rightarrow 0} \frac{\left.T\right|_{\params_0 + h} - \left.T\right|_{\params_0}}{h}.
	\end{equation}
	Both sensitivity measures are equivalent up to a change of sign and a scaling factor. The following relationship holds
	\begin{equation}
		S_T / T = - S_\omega / \omega.
	\end{equation}
\end{remark}

\subsection{Sensitivity analysis of a phase response curve}

Given a nominal parameter value $\params_0$, the sensitivity of the (input) infinitesimal phase response curve is the $2\pi$-periodic function $Z_{\iPRCu}:\Sbb^1\rightarrow\Rbb$ defined as
\begin{equation}
	Z_{\iPRCu}(\cdot) \eqdef \left.\frac{d\iPRCu}{d\params}(\cdot)\right|_{\params_0} = \lim_{h \rightarrow 0} \frac{\left.\iPRCu(\cdot)\right|_{\params_0 + h} - \left.\iPRCu(\cdot)\right|_{\params_0}}{h}.
\end{equation}
From~\eqref{eq:q_u}, we have, taking derivatives with respect to $\params$,
\begin{equation}
		Z_{\iPRCu}(\cdot) = \left\langle Z_{\iPRCx}(\cdot) , g(\tilde{x}^\gamma(\cdot),\params_0) \right\rangle + \left\langle \iPRCx(\cdot) , \frac{\partial g}{\partial x}(\tilde{x}^\gamma(\cdot),\params_0) Z_{\tilde{x}}(\cdot) + \frac{\partial g}{\partial \params}(\tilde{x}^\gamma(\cdot),\params_0)\right\rangle 
\end{equation}
where the $2\pi$-periodic function $Z_{\iPRCx}:\Sbb^1\rightarrow\Rbb^{n}$ is the sensitivity of the (state) infinitesimal phase response curve defined as
\begin{equation}
	Z_{\iPRCx}(\cdot) \eqdef \left.\frac{d\iPRCx}{d\params}(\cdot)\right|_{\params_0} = \lim_{h \rightarrow 0} \frac{\left.\iPRCx(\cdot)\right|_{\params_0 + h} - \left.\iPRCx(\cdot)\right|_{\params_0}}{h}.
\end{equation}
From~\eqref{eq:bvp_q}, we have, taking derivatives with respect to $\params$,
\begin{subequations} \label{eq:bvp_Zq}
	\begin{align}
		Z_{\iPRCx}'(\theta) + \frac{1}{\omega}\tr{A(\theta)} Z_{\iPRCx}(\theta) + \frac{1}{\omega} \tr{C(\theta)} \iPRCx(\theta) & = 0 \\
		Z_{\iPRCx}(2\pi) - Z_{\iPRCx}(0) & = 0 \\
		\left\langle Z_{\iPRCx}(\theta) , \tilde{v}(\theta) \right\rangle + \left\langle \iPRCx(\theta) , Z_{\tilde{v}}(\theta) \right\rangle  - S_\omega & = 0 
	\end{align}
\end{subequations}
where elements of the matrix $C(\cdot)$ are given by
\begin{align}
	C_{ij}(\cdot) & \eqdef \sum_{k=1}^{n} \frac{\partial^2f_i}{\partial x_j \partial x_k}(\tilde{x}^\gamma(\cdot),\params_0) (Z_{x})_k(\cdot) \nonumber \\
	& \qquad+ \frac{\partial^2 f_i}{\partial x_j \partial \params}(\tilde{x}^\gamma(\cdot),\params_0) - \frac{1}{\omega}\frac{\partial f_i}{\partial x_j}(\tilde{x}^\gamma(\cdot),\params_0)S_\omega ,
\end{align}
and where the $2\pi$-periodic function $Z_{\tilde{v}}:\Sbb^1\rightarrow\Rbb^{n}$ is the sensitivity of the vector field evaluated along the periodic orbit defined as
\begin{align}
	Z_{\tilde{v}}(\cdot) & \eqdef \left.\frac{d\tilde{v}}{d\params}(\cdot)\right|_{\params_0} =  \lim_{h \rightarrow 0} \frac{\left.\tilde{v}(\cdot)\right|_{\params_0 + h} - \left.\tilde{v}(\cdot)\right|_{\params_0}}{h} . 
\end{align}
Given the explicit dependence of the $2\pi$-periodic vector field in $\params$
\begin{equation}
	\tilde{v}(\cdot) = f(\left.\tilde{x}^\gamma(\cdot)\right|_{\params},\params),
\end{equation}
we have, taking derivatives with respect to $\params$, 
\begin{equation}
	Z_{\tilde{v}}(\cdot) = \frac{\partial f}{\partial x}(\tilde{x}^\gamma(\cdot),\params_0) Z_{\tilde{x}}(\cdot) + \frac{\partial f}{\partial \params}(\tilde{x}^\gamma(\cdot),\params_0).
\end{equation}

\subsection{Sensitivity analysis of the $1\mathord{:}1$~phase-locking} \label{sec:phase-lock-sens}

Given a nominal parameter value $\params_0$, the sensitivity of the phase difference $\chi^*$ is the scalar $S_{\chi^*}\in\Rbb$ defined as
\begin{equation}
	S_{\chi^*} \eqdef \lim_{h \rightarrow 0} \frac{\left.\chi^*\right|_{\params_0 + h} - \left.\chi^*\right|_{\params_0}}{h} .
\end{equation}
From $\chifield(\chi^*) = 0$, we have, taking derivatives of with respect to $\params$ and using \eqref{eq:phase_diff_dyn},
\begin{align}
	S_{\chi^*} & = -\left[\left.\chifield'\left(\left.\chi^*\right|_{\params_0}\right)\right|_{\params_0}\right]^{-1} \times \left[S_{\chifield}\left(\left.\chi^*\right|_{\params_0}\right)\right] \\
	& = - \left[\left.\Gamma'\left(\left.\chi^*\right|_{\params_0}\right)\right|_{\params_0}\right]^{-1} \times
	 \left[S_\Delta + S_{\Gamma}\left(\left.\chi^*\right|_{\params_0}\right)\right] 
\end{align}
where 
	$S_\Delta \eqdef\lim_{h\rightarrow0}[\Delta|_{\params_0+h} - \Delta|_{\params_0}]/h$ and 	
	$S_\Gamma(\cdot) \eqdef\lim_{h\rightarrow0}[\left.\Gamma(\cdot)\right|_{\params_0+h} - \left.\Gamma(\cdot)\right|_{\params_0}]/h$. 
Considering that $\omega|_\params = \Omega + \Delta|_\params$ is the sum of a parameter independent term~$\Omega$ and a parameter dependent term $\Delta$, we have that $S_\omega = S_\Delta$. In addition, from~\eqref{eq:coupl_fun}, we have, taking derivatives with respect to $\params$,
\begin{equation}
	S_{\Gamma}(\cdot) = \lim_{\tilde{T}\rightarrow+\infty}\frac{1}{\tilde{T}} \int_0^{\tilde{T}} S_{\iPRCu}(\Omega t + \cdot) \tilde{h}_u(\Omega_u t ) dt .
\end{equation}
The sensitivity of the phase difference has thus two distinct contributions: 
\begin{equation}
	S_{\chi^*} = S_{\chi^*|\omega} + S_{\chi^*|\Gamma}
\end{equation}
where $S_{\chi^*|\omega}\eqdef- [\Gamma'(\chi^*|_{\params_0})|_{\params_0}]^{-1} \times S_\omega$ denotes the contribution of the angular frequency sensitivity and $S_{\chi^*|\Gamma}\eqdef- [\Gamma'(\chi^*|_{\params_0})|_{\params_0}]^{-1} \times S_{\Gamma}(\chi^*|_{\params_0})$ denotes the contribution of the coupling function sensitivity at $\chi^*$, the latter being closely related to the iPRC.

\subsection{Numerics of sensitivity analysis} \label{sec:numerics}

Numerical algorithms to solve boundary value problems \eqref{eq:bvp_Zx} and \eqref{eq:bvp_Zq} are reviewed in \cite[Appendix]{Sacre:2012uz}. We stress that existing algorithms that compute periodic orbits and iPRCs are easily adapted to compute their sensitivity curves, essentially at the same numerical cost. All numerical tests in Section~\ref{sec:illustations} have been obtained with a MATLAB numerical code available from the authors. 

The proposed approach is systematic and computationally tractable but it only provides a \emph{local} sensitivity analysis in the parameter space, around a nominal set of parameter values. It complements more \emph{global}---but less tractable---tools such as bifurcation analysis or parameter space exploration. Studying the bifurcation diagram associated with a given parameter~\cite{Leloup:2004co} or using sampling methods in the full parameter space~\cite{Stelling:2004do,Hafner:2009eb} are classical ways to assess the robustness of an oscillator: the parameter range over which the oscillation exists is a (nonlocal) indicator of the sensitivity of the oscillator to the parameters. The limitation of those approaches is that they are univariate (only one direction of the parameter space is explored in a particular bifurcation diagram) and that the exploration of the parameter space rapidly becomes formidable as the number of parameters grows.


\section{Scalar robustness measures for oscillators} 
\label{sec:scalar_robustness_measures_for_oscillators}

Testing the robustness of a model against parameter variations is a basic system-theoretic question. In a number of situations, the very purpose of modeling is to identify those parameters that influence a given system property. In the literature, robustness analysis of circadian rhythms mostly studies the zero-input steady-state behavior such as the period or the amplitude of oscillations~\cite{Gonze:2002im,Stelling:2004do,Wilkins:2007is} and (empirical) phase-based performance measures~\cite{Bagheri:2007bo,Gunawan:2007jv,Hafner:2010us,Pfeuty:2011em}. In this section, we propose scalar robustness measures to quantify the sensitivity of  the angular frequency,  the infinitesimal phase response curve, and  the $1\mathord{:}1$~phase-locking to parameters. 

\subsubsection*{Robustness measure of the angular frequency}
The angular frequency~$\omega$ is a positive scalar number. The sensitivity of~$\omega$ with respect to the parameter~$\params$ is thus also a scalar number~$S_\omega$, leading to a scalar robustness measure~$R_\omega$ defined as
\begin{equation}
	R_{\omega} \eqdef \left| S_{\omega} \right|
\end{equation}
where $|\cdot|$ denotes the real absolute value function.

\subsubsection*{Robustness measure of the infinitesimal phase response curve}

In contrast, the iPRC $\iPRCu:\Sbb^1\rightarrow\Rbb$ belongs to an infinite-dimensional space~$\mathcal{Q}$. The sensitivity of $\iPRCu$ with respect to the parameter $\params$ is thus a vector~$S_{\iPRCu}$ which belongs to the tangent space $\tgspace{\iPRCu}{\mathcal{Q}}$ at $\iPRCu$. A scalar robustness measure~$R_{\iPRCu}$ is defined as
\begin{equation}
	R_{\iPRCu} \eqdef \left\| S_{\iPRCu} \right\|_{\iPRCu} = \sqrt{g_{\iPRCu}\left(S_{\iPRCu},S_{\iPRCu}\right)}
\end{equation}
where $\left\| \cdot \right\|_{\iPRCu}$ denotes the norm induced by a Riemannian metric $g_{\iPRCu}\left(\cdot,\cdot\right)$ at $\iPRCu$. 
In this chapter, we use the simplest metric for signals in $\mathcal{L}_2(\Sbb^1,\Rbb)$, that is the standard inner product,
\begin{equation}
	g_{\iPRCu}( \xi_{\iPRCu},\zeta_{\iPRCu} ) \eqdef \left\langle \xi_{\iPRCu},\zeta_{\iPRCu} \right\rangle = \int_{\Sbb^1} \xi_{\iPRCu}(\theta) \conj{\zeta_{\iPRCu}(\theta)} d\theta .
\end{equation}
In our previous work \cite{Sacre:2012uz}, we proposed further metrics which capture equivalence properties in the space of phase response curves. This is motivated by the fact that, in many applications, it is not meaningful to distinguish among PRCs that are related by a scaling factor and/or a phase shift. 

\subsubsection*{Robustness measure of the $1\mathord{:}1$~phase-locking}

The stable phase difference $\chi^*$ is a scalar phase on the unit circle $\Sbb^1$. The sensitivity of $\chi^*$ with respect to the parameter $\params$ is a scalar number $S_{\chi^*}$, leading to a scalar robustness measure $R_{\chi^*}$ defined as
\begin{equation}
	R_{\chi^*} \eqdef \left| S_{\chi^*} \right|.
\end{equation}

\subsubsection*{Normalized robustness measures}

When analyzing a model with several parameters ($\params \in \Params \subseteq \Rbb^q$), all robustness measures $R_\star$ (where $\star$ stands for any characteristic of the oscillator) are $q$-dimensional vectors. Each element of those vectors represents to the scalar robustness measure corresponding to each parameter.
The normalized robustness measure
\begin{equation}
	\overline{R}_{\star} = \frac{R_\star}{\left\|R_\star\right\|_{\infty}}
\end{equation}
has all its components in the unit interval~$[0,1]$. This normalized measure allows to rank model parameters according to their relative ability to influence the characteristic $\star$.

\begin{remark} \label{rem:period_norm_rob}
	Note that $\overline{R}_{T} = \overline{R}_{\omega}$.
\end{remark}


\section{Application to a model of circadian rhythms} 
\label{sec:illustations}

\begin{figure}
	\centering
	\includegraphics[width=\textwidth]{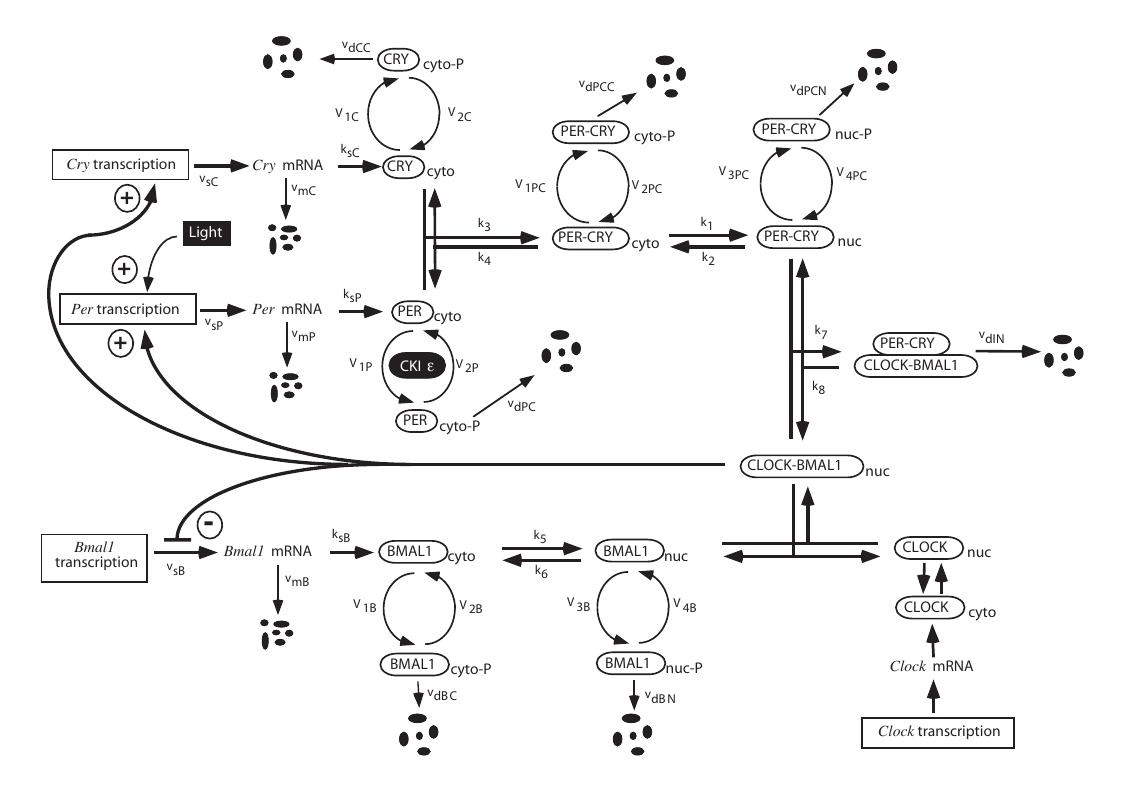}
	\caption{The Leloup-Goldbeter model accounts for several regulatory processes identified in circadian rhythms of mammals. Reproduction of a figure from \protect\cite{Leloup:2003cp}.}
	\label{fig:goldbeter-model}
\end{figure}

We illustrate our sensitivity analysis on the genetic oscillator model of Leloup and Goldbeter (see \myfigurename~\ref{fig:goldbeter-model}). This model accounts for several regulatory processes identified in circadian rhythms of mammals. 
A negative autoregulatory feedback loop established by the \textit{per} (period) and \textit{cry} (cryptochrome) genes is at the heart of the circadian oscillator.  The PER and CRY proteins form a complex PER--CRY that indirectly represses the activation of the \textit{Per} and \textit{Cry} genes. The PER--CRY complexes exert their repressive effect by binding to a complex of two proteins CLOCK--BMAL1. This latter, formed by the products of \textit{Clock} and \textit{Bmal1} genes, activates \textit{Per} and \textit{Cry} transcription. 
In addition to this negative autoregulation, an (indirect) positive regulatory feedback loop is also involved. Indeed, the \textit{Bmal1} expression is subjected to negative autoregulation by CLOCK--BMAL1, through the product of the \textit{Rev-Erb$\alpha$} gene. The complex PER--CRY enhances \textit{Bmal1} expression in an indirect manner by binding to CLOCK--BMAL1, and thereby reducing the transcription of the \textit{Rev-Erb$\alpha$} gene.
Finally, environmental periodic cycles associated with earth's rotation are mediated through light--dark cycles. Light acts on the system by inducing the expression of the \textit{Per} gene. 

The detailed computational model of Leloup and Goldbeter possesses $16$~state variables and $52$~parameters. State-space model equations and nominal parameter values are available in \cite[Supporting Text]{Leloup:2003cp}. The effect of light is incorporated through periodic square-wave variations in the maximal rate of \textit{Per} expression (\ie{} the value of the parameter $v_{\text{sP}}$ goes from a constant low value during dark phase to a constant  high value during light phase). Parameters values remain to be determined experimentally and have been chosen semiarbitrarily in physiological ranges in order to satisfy experimental observations.

Each parameter of the model describes a single regulatory mechanism such as transcription and translation control of mRNAs, degradation of mRNAs or proteins, transport reaction, and phosphorylation/dephosphorylation of proteins. 
The analysis of single-parameter sensitivities reveals thus the importance of individual regulatory processes on the function of the oscillator. 

However, in order to enlighten the potential role of circuits rather than single-parameter properties, we grouped model parameters according to the mRNA loop to which they belong. Each group of parameters is associated with a different color: \textit{Per}-loop in blue, \textit{Cry}-loop in red, and \textit{Bmal1}-loop in green. In addition, we gathered parameters associated with interlocked loops in a last group represented in gray.

\emph{In the following, we consider sensitivities to relative variations of parameters. We write without distinction about period sensitivities and angular frequency sensitivities due to their direct proportional relationship (see remarks~\ref{rem:period_sens} and \ref{rem:period_norm_rob}).}

\subsubsection*{Sensitivity analysis of the period and the phase response curve}

The period and the PRC are two intrinsic characteristics of the circadian oscillator with physiological significance. We use the sensitivity analysis of the period and the PRC to measure the influence of regulatory processes on tuning the period and shaping the PRC.

\begin{figure}
	\centering
	\includegraphics{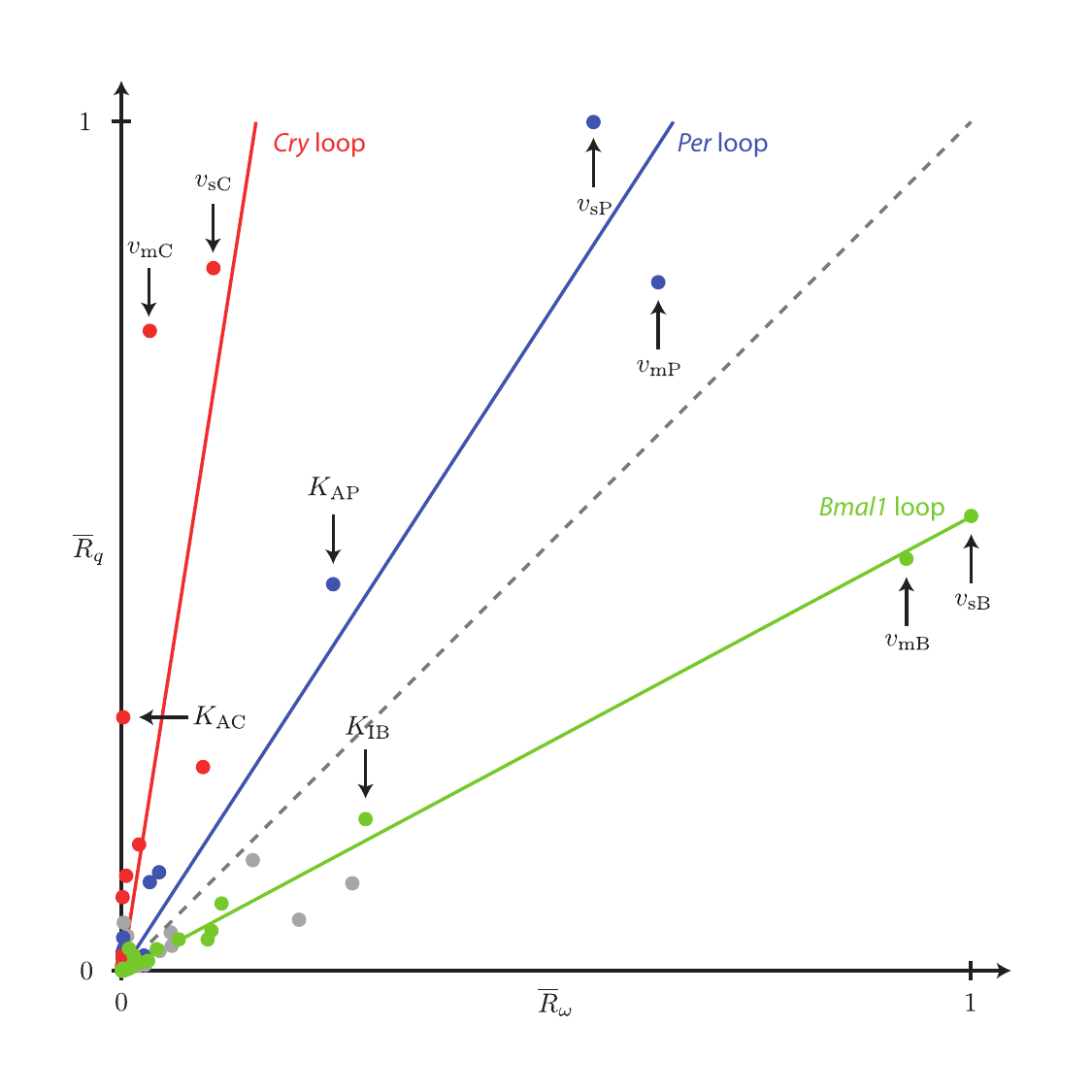}
	\caption{Normalized robustness measures $\overline{R}_\omega$ (angular frequency) and $\overline{R}_{\iPRCu}$ (iPRC) reveal the distinct sensitivity of three distinct genetic circuits (\emph{Cry}, \emph{Per}, and \emph{Bmal1}). 
	Each point is associated to a particular parameter.
	The three lines are regression over the parameters of the three gene loops.
	The dashed bisector indicates the positions at which both measures of robustness are identical.
	Only parameters associated with the \emph{Cry}-loop exhibit low angular frequency and high iPRC sensitivities. 
	The color code corresponds to different subsets of parameters associated to different loops (see the text for details).}
	\label{fig:goldbeter-results1}
\end{figure}

A two-dimensional $(\overline{R}_\omega, \overline{R}_{\iPRCu})$ scatter plot in which each point corresponds to a parameter of the model reveals the shape and  strength of the relationship between both normalized robustness measures $\overline{R}_\omega$ (angular frequency or, equivalently, period) and $\overline{R}_{\iPRCu}$ (PRC). It enables to identify which characteristic is primarily affected by perturbations in individual parameters: parameters corresponding to points situated below the dashed bisector influence mostly the period; those above the dashed bisector influence mostly the PRC (see \myfigurename~\ref{fig:goldbeter-results1}).

At a coarse level of analysis, the scatter plot reveals that most parameters exhibit both low period and PRC sensitivities (most points are close to the origin); only few parameters display a medium or high sensitivity either to period or to PRC.

At a finer level of analysis, the scatter plot  reveals that the parameters associated with each of the three mRNA loops have distinct sensitivities:
\begin{itemize}
	\item the \emph{Bmal1}-loop parameters are associated with a high period sensitivity and a medium PRC sensitivity (regression line below the bisector);
	\item the \emph{Per}-loop parameters are associated with a medium period sensitivity and a high PRC sensitivity (regression line above the bisector);
	\item the \emph{Cry}-loop parameters are associated with a low period sensitivity and a high PRC sensitivity (regression line above the bisector, close to the vertical axis).
\end{itemize}
In each feedback loop, the three more sensitive parameters represent the three same biological functions: the maximum rates of mRNA synthesis ($v_\text{sB}$, $v_\text{sP}$, and $v_\text{sC}$), the maximum rate of mRNA degradation ($v_\text{mB}$, $v_\text{mP}$, and $v_\text{mC}$), and the inhibition~(I) or activation~(A) constants for the repression or enhancement of mRNA expression by BMAL1 ($K_\text{IB}$, $K_\text{AP}$, and~$K_\text{AC}$).

The small number of highly sensitive parameters is in agreement with the robust nature of the circadian clock and the concentration of fragilities in some specific locations of the architecture~\cite{Stelling:2004do}. Our analysis suggests that the transcriptional and translational control of mRNA (\ie{} the control of both biological steps required to synthesize a protein) has to be regulated by specific mechanisms (not included in the model) in order to avoid failures in the clock function.
While the topology of \textit{Per}- and \textit{Cry}-loops are identical, the asymmetry introduced by the choice of parameter values leads to different sensitivity for those loops. Both loops have a similar high sensitivity of the PRC (while the light acts only on the maximum rate of Per mRNA synthesis) but a different sensitivity of the period, the \textit{Per}-loop being more sensitive than the \textit{Cry}-loop.
The high sensitivity of the period for parameters associated with the \textit{Bmal1}-loop has also being identified in~\cite{Leloup:2004co}. However, this last prediction of the model (high sensitivity of the period to \textit{Bmal1}-loop) is not in agreement with  experimental observations in~\cite{Bunger:2000cs,vonGall:2003br}. This observation may encourage the biologist and the modeler to design of new experiments to enlighten biological mechanisms responsible for this discrepancy between the experiment and the model.

\subsubsection*{Sensitivity analysis of the entrainment}

Entrainment is an important characteristic of the circadian model. In Section~\ref{sec:phase-lock-sens}, we have seen that the entrainment sensitivity $S_{\chi^*}$ is mathematically given by the summation of two terms: a term $S_{\chi^*|\omega}$ proportional to the period sensitivity and a term $S_{\chi^*|\Gamma}$ proportional to the coupling function sensitivity at $\chi^*$. Those two terms correspond to two biologically distinct mechanisms by which the entrainment properties of the circadian clock can be regulated: a modification of the period or a modification of the coupling function (resulting from the modification of the iPRC or the input signal).

\begin{figure}
	\centering
	\includegraphics{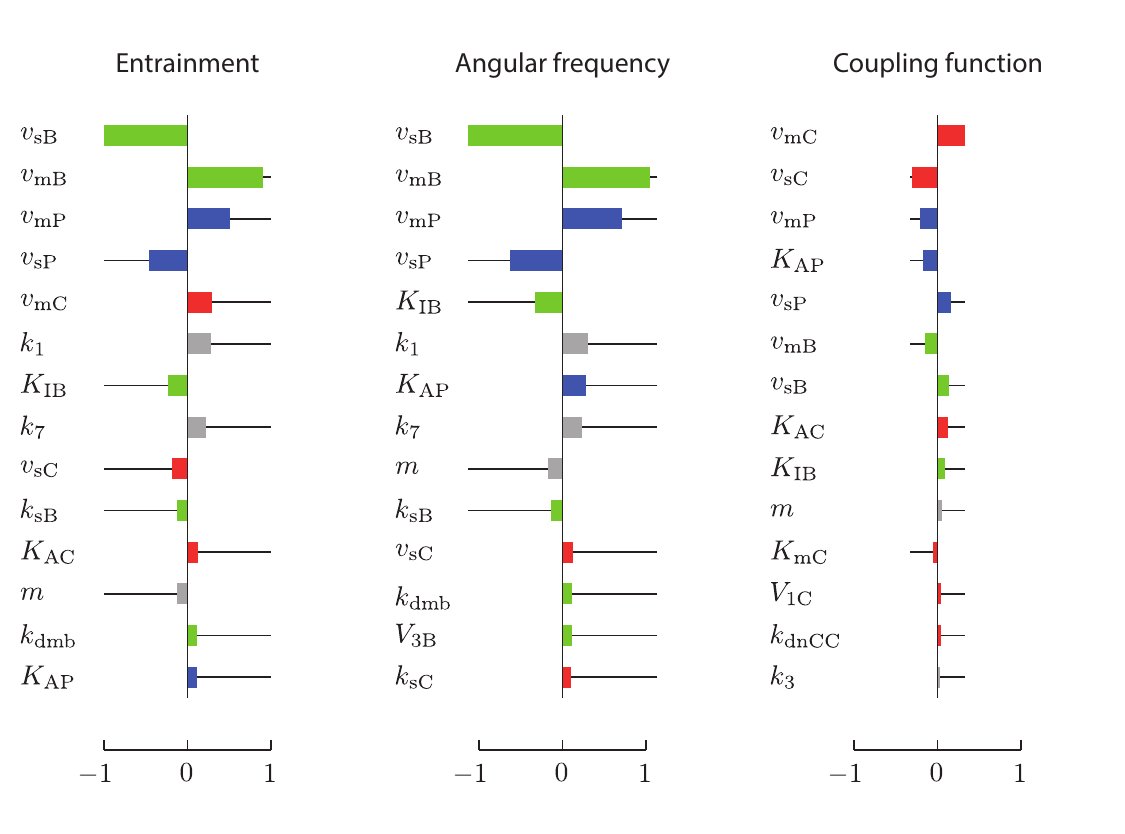}
	\caption{Normalized sensitivity measures $S_{\chi^*}/\left\|S_{\chi^*}\right\|_{\infty}$ (entrainment) are due to two contributions: $S_{\chi^*|\omega}/\left\|S_{\chi^*}\right\|_{\infty}$ (angular frequency) and $S_{\chi^*|\Gamma}/\left\|S_{\chi^*}\right\|_{\infty}$ (coupling function). Each (thick) horizontal bar corresponds to a sensitivity measure with respect to a particular parameter. The (thin) horizontal lines indicate (in absolute value) the maximal sensitivity (among all parameters) and may be useful to compare the sensitivity of a parameter to the maximal sensitivity.
	The color code corresponds to different subsets of parameters associated to different loops (see the text for details).}
	\label{fig:goldbeter-results2}
\end{figure}

Bar plots of $S_{\chi^*}/\left\|S_{\chi^*}\right\|_{\infty}$, $S_{\chi^*|\omega}/\left\|S_{\chi^*}\right\|_{\infty}$, and $S_{\chi^*|\Gamma}/\left\|S_{\chi^*}\right\|_{\infty}$ in which each bar corresponds to a parameter allows to identify the most sensitive parameters for entrainment and to quantify\footnote{The entrainment sensitivity and the  contributing terms are normalized by $\left\|S_{\chi^*}\right\|_{\infty}$ (the same maximal value of the entrainment sensitivity) such that the summation of normalized terms is equal to the normalized entrainment sensitivity.} the respective contribution of both mechanisms in the entrainment sensitivity (see \myfigurename~\ref{fig:goldbeter-results2}). For each bar plot, we sorted parameters by absolute magnitude and restricted the plot to the 14 parameters with the highest sensitivity measure (the number 14 results from our choice to keep the parameters with an entrainment sensitivity greater than 0.1). Those plots allow to identify the parameters which play an important role in the entrainment sensitivity. We note that the parameter orders for $S_{\chi^*}/\left\|S_{\chi^*}\right\|_{\infty}$ and $S_{\chi^*|\omega}/\left\|S_{\chi^*}\right\|_{\infty}$ are almost identical, except for parameters associated
with the \textit{Cry}-loop. Those parameters appear in the highest ones for $S_{\chi^*|\Gamma}/\left\|S_{\chi^*}\right\|_{\infty}$.

\begin{figure}
	\centering
	\includegraphics{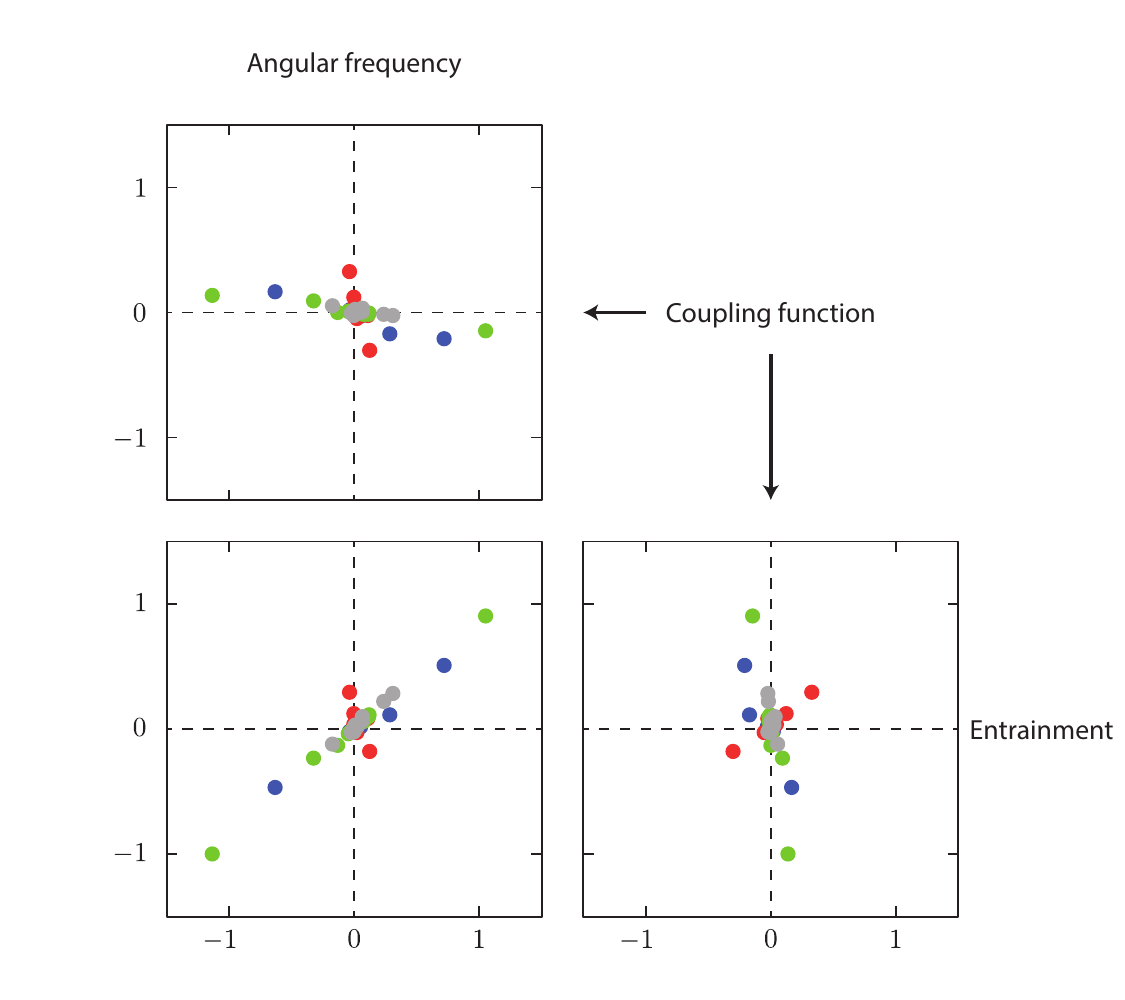}
	\caption{Normalized sensitivity measures $S_{\chi^*}/\left\|S_{\chi^*}\right\|_{\infty}$ (entrainment), $S_{\chi^*|\omega}/\left\|S_{\chi^*}\right\|_{\infty}$ (angular frequency), and $S_{\chi^*|\Gamma}/\left\|S_{\chi^*}\right\|_{\infty}$ (coupling function) exhibit particular correlation shapes. 
	The top graph represents the  $(S_{\chi^*|\omega}/\left\|S_{\chi^*}\right\|_{\infty},S_{\chi^*|\Gamma}/\left\|S_{\chi^*}\right\|_{\infty})$-plan; the bottom-left graph represents the $(S_{\chi^*|\omega}/\left\|S_{\chi^*}\right\|_{\infty},S_{\chi^*}/\left\|S_{\chi^*}\right\|_{\infty})$-plan; and the bottom-right graph represents the $(S_{\chi^*|\Gamma}/\left\|S_{\chi^*}\right\|_{\infty},S_{\chi^*}/\left\|S_{\chi^*}\right\|_{\infty})$-plan.
	Each point is associated to a particular parameter.
	The color code corresponds to different subsets of parameters associated to different loops (see the text for details).
	Those correlations support the competitive nature of both mechanisms (modification of the period or the coupling function) leading to the entrainment sensitivity.}
	\label{fig:goldbeter-results3}
\end{figure}

Figure~\ref{fig:goldbeter-results3} (top) reveals the competitive and complementary nature of both contributions to entrainment sensitivity. For most parameters, both contributions have opposite signs, that is, points are located in the second and fourth quadrants. In addition, both mechanisms are well decoupled such that, when one mechanism is active, the other is almost inactive (points are located close to the horizontal and vertical axes). 
Parameters associated with \textit{Cry}-loop seem to influence the entrainment sensitivity through a modification of the coupling function (points close to the vertical axis); others parameters associated with \textit{Per}-loop and \textit{Bmal1}-loop seem to influence the entrainment sensitivity through a modification of the period (points close to the horizontal axis). 

The different mechanisms leading to entrainment sensitivity  are also observed in both other scatter plots (see \myfigurename~\ref{fig:goldbeter-results3} bottom-left and -right). In those plots, parameters associated with points close to the bisector of the first and third quadrants influence the entrainment sensitivity through a modification of the period (bottom-left) or the coupling function (bottom-right), respectively. Again, only parameters associated with the \emph{Cry}-loop seem to affect the entrainment through a variation of the PRC.

Two of the parameters belonging to the \emph{Cry}-loop (with high coupling function and low period sensitivities) have been identified by numerical simulations as important for entrainment properties of the model without affecting the period: $K_{\text{AC}}$~in~\cite{Leloup:2003cp} and $v_{\text{mC}}$~in~\cite{Leloup:2004co}. Our approach supports the importance
of those two parameters and identifies the potential importance of a third
one~($v_{\text{sC}}$).

\bigskip

We stress that the sensitivity analysis in \cite{Leloup:2003cp,Leloup:2004co} is a \emph{global} approach that relies on exploring the parameter space through numerical simulations of the model to determine the system behavior under constant and periodic environmental conditions while varying one parameter at a time. In contrast, the proposed analysis is \emph{local} but systematic and computationally tractable.
In the particular model studied here and in \cite{Leloup:2004co}, the predictions of the (local) sensitivity analysis match the predictions of the (nonlocal) analysis.

\begin{figure}
	\centering
	\includegraphics[width=11.43cm]{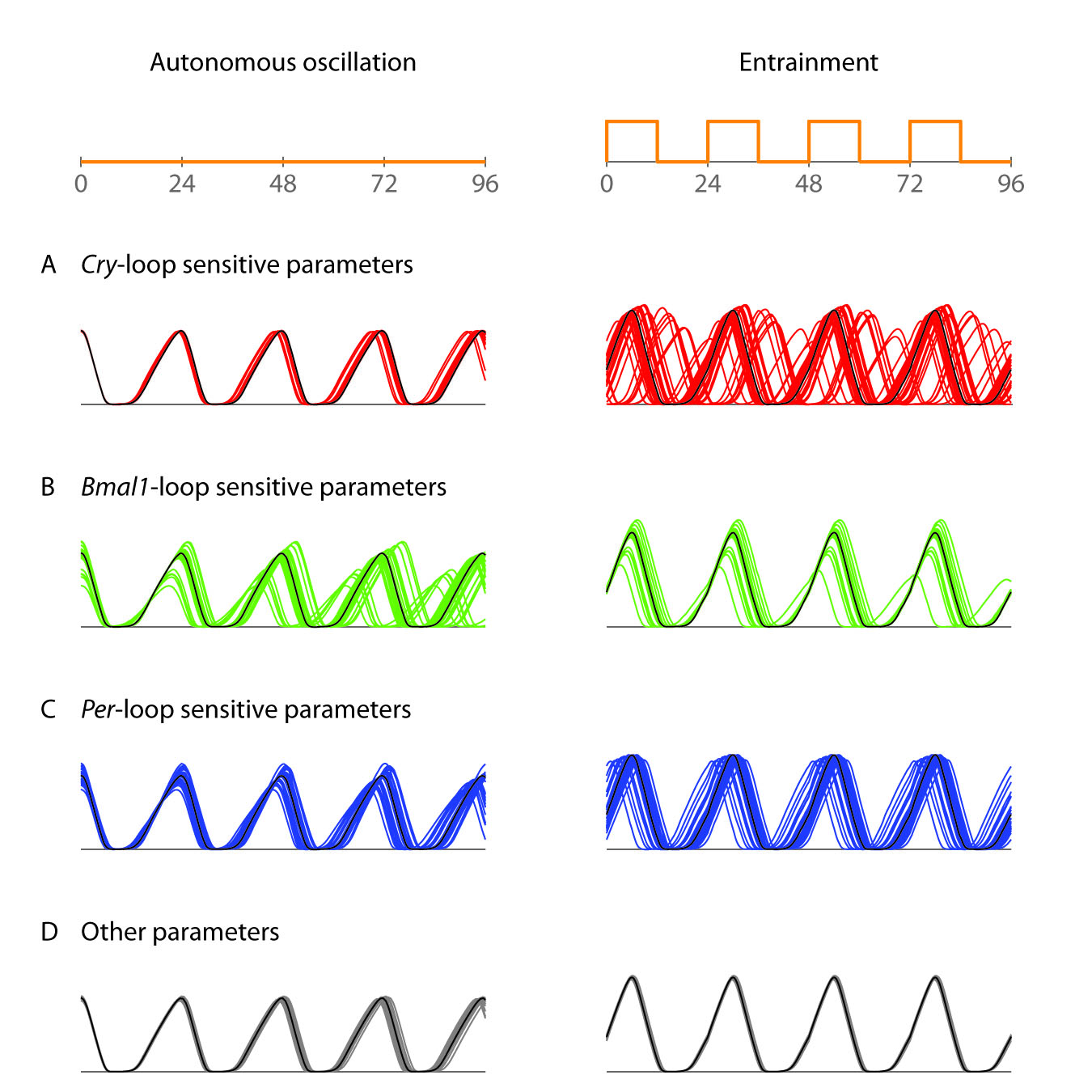}
	\caption{%
	Steady-state behaviors for the nominal model and different finite (nonlocal) parameter perturbations are illustrated by time-plots of the state variable $M_P$ under constant environmental conditions (autonomous oscillation, \emph{left}) and periodic environmental conditions (entrainment, \emph{right}).
	 Each panel (or row) corresponds to the perturbation of a different group of parameters, the black time-plot corresponding to system behaviors for nominal parameter values. Perturbations are randomly taken in a range of $\pm 10\%$ around the nominal parameter value (for one parameter at a~time).
	\textbf{A.}~Perturbations of three most sensitive parameters of \emph{Cry}-loop ($v_{\text{sC}}$, $v_{\text{mC}}$, and $K_{\text{AC}}$) lead to small variations of the autonomous period and (not structured) large variations of the phase-locking.
	\textbf{B.}~Perturbations of three most sensitive parameters of \emph{Bmal1}-loop ($v_{\text{sB}}$, $v_{\text{mB}}$, and $K_{\text{IB}}$) lead to larger variations of the autonomous period and medium variations of the phase-locking.
	\textbf{C.}~Perturbations of three most sensitive parameters of \emph{Per}-loop ($v_{\text{sP}}$, $v_{\text{mP}}$, and $K_{\text{AP}}$) exhibit an intermediate behavior between the situations A and B.
	\textbf{D.}~Perturbations of parameters of interlocked loops lead to small variations of the autonomous period and the phase-locking.
	}
	\label{fig:global-sensitivity}
\end{figure}

To evaluate the nonlocal nature of our local predictions, we plot in \myfigurename~\ref{fig:global-sensitivity} the time behavior of solutions for different finite (nonlocal) parameter changes. The left plots illustrate the autonomous oscillation of the isolated oscillator whereas the right plots illustrate the steady-state solution entrained by a periodic light input.
Parameter perturbations are randomly taken in a range of $\pm 10\%$ around the nominal parameter value. Each panel corresponds to the perturbation of a different group of parameters (the black time-plot corresponds to the nominal system behaviors for nominal parameter values).
\begin{enumerate}[A.]
	\item Perturbations of three most sensitive parameters of \emph{Cry}-loop ($v_{\text{sC}}$, $v_{\text{mC}}$, and $K_{\text{AC}}$) lead to small variations (mostly shortening) of the autonomous period and (not structured) large variations of the phase-locking. This observation is consistent with the low sensitivity of the period and the high sensitivity of the PRC. 
	\item Perturbations of three most  sensitive parameters of \emph{Bmal1}-loop ($v_{\text{sB}}$, $v_{\text{mB}}$, and $K_{\text{IB}}$) lead to medium variations of the autonomous period and medium variations of the phase-locking. The variations of the phase-locking exhibit the same structure as variations of the period, suggesting that the change in period is responsible for the change of phase-locking for those parameters. This observation is consistent with the high sensitivity of the period and the medium sensitivity of the PRC. 	
	\item Perturbations of three most  sensitive parameters of \emph{Per}-loop ($v_{\text{sP}}$, $v_{\text{mP}}$, and $K_{\text{AP}}$) exhibit an intermediate behavior between the situations A and B.
	\item Perturbations of parameters of interlocked loops lead to small variations of the autonomous period and the phase-locking, which is consistent with their low sensitivity.
\end{enumerate}
Those (nonlocal) observations are thus well predicted by the classification of parameters suggested by the (local) sensitivity analysis (see \myfigurename~\ref{fig:goldbeter-results1}).

\section{Conclusion} 
\label{sec:conclusion}

This chapter proposes (local) sensitivity tools to analyze oscillator models as open dynamical systems. 
We showed that, under the weak perturbation assumption, state-space models can be reduced to phase models characterized by their angular frequency and their phase response curve. Those phase models are then useful to study the entrainment (or phase-locking) to a periodic input.
We then introduced the sensitivity analysis for oscillators and their phase-locking behavior.

The application of this approach to a detailed computational model of circadian rhythms provides physiologically relevant predictions. It enlightens the distinct role of different circuits in the robustness of entrainment and it selects 3 out of 52 parameters as parameters that strongly affect the phase response curve while barely affecting the period. The importance of two of these parameters was previously identified in the literature through simulations of the model.


\section{Lessons learnt} 
\label{sec:lessons_learnt}

Sensitivity analysis is a classical and fundamental tool to evaluate the role of a given parameter in a given system characteristic. Because the phase response curve is a fundamental input--output characteristic of oscillators, we developed a sensitivity analysis for oscillator models in the space of phase response curves. The proposed tool can be applied to high-dimensional  oscillator models without facing the curse of dimensionality obstacle associated with numerical exploration of the parameter space. Application of this tool to a state-of-the-art model of circadian rhythms suggests that it can be useful and instrumental to biological investigations.


\bibliographystyle{unsrt}
\bibliography{bib/bibfile}

\begin{thebibliography}{10}

\bibitem{Pittendrigh:1981wa}
Colin~S Pittendrigh.
\newblock {Circadian systems: entrainment}.
\newblock In {\em Biological Rhythms}, pages 95--124. Plenum Press, New York,
  NY, 1981.

\bibitem{Hastings:2000ht}
Michael~H Hastings.
\newblock {Circadian clockwork: two loops are better than one}.
\newblock {\em Nat. Rev. Neurosci.}, 1(2):143--146, November 2000.

\bibitem{Stelling:2004do}
J{\"o}rg Stelling, Ernst~Dieter Gilles, and Francis~J Doyle~III.
\newblock {Robustness properties of circadian clock architectures}.
\newblock {\em Proc. Natl. Acad. Sci. U.S.A.}, 101(36):13210--13215, September
  2004.

\bibitem{Novak:2008eg}
B{\'e}la Nov{\'a}k and John~J Tyson.
\newblock {Design principles of biochemical oscillators}.
\newblock {\em Nat. Rev. Mol. Cell Biol.}, 9(12):981--991, December 2008.

\bibitem{Winfree:1967vf}
Arthur~T Winfree.
\newblock {Biological rhythms and the behavior of populations of coupled
  oscillators}.
\newblock {\em J. Theor. Biol.}, 16(1):15--42, July 1967.

\bibitem{Winfree:1980ue}
Arthur~T Winfree.
\newblock {\em {The Geometry of Biological Time}}, volume~8 of {\em
  Biomathematics}.
\newblock Springer-Verlag, New York, NY, 1st edition, 1980.

\bibitem{Leloup:2003cp}
Jean-Christophe Leloup and Albert Goldbeter.
\newblock {Toward a detailed computational model for the mammalian circadian
  clock}.
\newblock {\em Proc. Natl. Acad. Sci. U.S.A.}, 100(12):7051--7056, June 2003.

\bibitem{Leloup:2004co}
Jean-Christophe Leloup and Albert Goldbeter.
\newblock {Modeling the mammalian circadian clock: sensitivity analysis and
  multiplicity of oscillatory mechanisms}.
\newblock {\em J. Theor. Biol.}, 230(4):541--562, October 2004.

\bibitem{Hafner:2009eb}
Marc Hafner, Heinz Koeppl, Martin Hasler, and Andreas Wagner.
\newblock {`Glocal' robustness analysis and model discrimination for circadian
  oscillators}.
\newblock {\em PLoS Comput. Biol.}, 5(10):e1000534, October 2009.

\bibitem{Sepulchre:2006vk}
Rodolphe Sepulchre.
\newblock {Oscillators as systems and synchrony as a design principle}.
\newblock In {\em Current Trends in Nonlinear Systems and Control: In Honor of
  Petar Kokotovi{\'c} and Turi Nicosia}, pages 123--141. Birkh{\"a}user,
  Boston, MA, 2006.

\bibitem{Farkas:1994uq}
Miklos Farkas.
\newblock {\em {Periodic Motions}}, volume 104 of {\em Applied Mathematical
  Sciences}.
\newblock Springer-Verlag, New York, NY, 1994.

\bibitem{Khalil:2002wj}
Hassan~K Khalil.
\newblock {\em {Nonlinear Systems}}.
\newblock Prentice Hall, Upper Saddle River, NJ, 3rd edition, 2002.

\bibitem{Kuramoto:1984wo}
Yoshiki Kuramoto.
\newblock {\em {Chemical Oscillations, Waves, and Turbulence}}, volume~19 of
  {\em Springer Series in Synergetics}.
\newblock Springer-Verlag, Berlin/Heidelberg, Germany, 1st edition, 1984.

\bibitem{Hoppensteadt:1997tp}
Frank~C Hoppensteadt and Eugene~M Izhikevich.
\newblock {\em {Weakly Connected Neural Networks}}, volume 126 of {\em Applied
  Mathematical Sciences}.
\newblock Springer-Verlag, New York, NY, 1997.

\bibitem{Ascher:1988ty}
Uri~M Ascher, Robert M~M Mattheij, and Robert~D Russell.
\newblock {\em {Numerical Solution of Boundary Value Problems for Ordinary
  Differential Equations}}.
\newblock Prentice Hall, Englewood Cliffs, NJ, February 1988.

\bibitem{Seydel:2010tw}
R{\"u}diger Seydel.
\newblock {\em {Practical Bifurcation and Stability Analysis}}, volume~5 of
  {\em Interdisciplinary Applied Mathematics}.
\newblock Springer, New York, NY, 3rd edition, 2010.

\bibitem{Sacre:2012uz}
Pierre Sacr{\'e} and Rodolphe Sepulchre.
\newblock {System analysis of oscillator models in the space of phase response
  curves}.
\newblock {\em arXiv}, math.DS, June 2012.

\bibitem{Malkin:1949to}
Ioel'~Gil'evich Malkin.
\newblock {\em {The Methods of Lyapunov and Poincare in the Theory of Nonlinear
  Oscillations}}.
\newblock Gostexizdat, Moscow, Russia, 1949.

\bibitem{Malkin:1956wq}
Ioel'~Gil'evich Malkin.
\newblock {\em {Some Problems in the Theory of Nonlinear Oscillations}}.
\newblock Gostexizdat, Moscow, Russia, 1956.

\bibitem{Neu:1979cj}
John~C Neu.
\newblock {Coupled chemical oscillators}.
\newblock {\em SIAM J. Appl. Math.}, 37(2):307--315, 1979.

\bibitem{Ermentrout:1984jb}
G~Bard Ermentrout and Nancy Kopell.
\newblock {Frequency plateaus in a chain of weakly coupled oscillators. I}.
\newblock {\em SIAM J. Math. Anal.}, 15(2):215--237, 1984.

\bibitem{Kramer:1984jk}
Mark~A Kramer, Herschel Rabitz, and Joseph~M Calo.
\newblock {Sensitivity analysis of oscillatory systems}.
\newblock {\em Appl. Math. Modelling}, 8(5):328--340, October 1984.

\bibitem{Rosenwasser:1999tw}
Efim Rosenwasser and Rafael Yusupov.
\newblock {\em {Sensitivity of Automatic Control Systems}}.
\newblock CRC Press, Boca Raton, FL, 1999.

\bibitem{Ingalls:2004wr}
Brian~P Ingalls.
\newblock {Autonomously oscillating biochemical systems: parametric sensitivity
  of extrema and period}.
\newblock {\em Syst. Biol.}, 1(1):62--70, June 2004.

\bibitem{Wilkins:2009kq}
Anna~Katharina Wilkins, Bruce Tidor, Jacob White, and Paul~I Barton.
\newblock {Sensitivity analysis for oscillating dynamical systems}.
\newblock {\em SIAM J. Sci. Comput.}, 31(4):2706--2732, 2009.

\bibitem{Gonze:2002im}
Didier Gonze, Jos{\'e} Halloy, and Albert Goldbeter.
\newblock {Robustness of circadian rhythms with respect to molecular noise}.
\newblock {\em Proc. Natl. Acad. Sci. U.S.A.}, 99(2):673--678, January 2002.

\bibitem{Wilkins:2007is}
Anna~Katharina Wilkins, Paul~I Barton, and Bruce Tidor.
\newblock {The Per2 negative feedback loop sets the period in the mammalian
  circadian clock mechanism}.
\newblock {\em PLoS Comput. Biol.}, 3(12):e242, December 2007.

\bibitem{Bagheri:2007bo}
Neda Bagheri, J{\"o}rg Stelling, and Francis~J Doyle~III.
\newblock {Quantitative performance metrics for robustness in circadian
  rhythms}.
\newblock {\em Bioinformatics}, 23(3):358--364, February 2007.

\bibitem{Gunawan:2007jv}
Rudiyanto Gunawan and Francis~J Doyle~III.
\newblock {Phase sensitivity analysis of circadian rhythm entrainment}.
\newblock {\em J. Biol. Rhythms}, 22(2):180--194, April 2007.

\bibitem{Hafner:2010us}
Marc Hafner, Pierre Sacr{\'e}, Laura Symul, Rodolphe Sepulchre, and Heinz
  Koeppl.
\newblock {Multiple feedback loops in circadian cycles: robustness and
  entrainment as selection criteria}.
\newblock In {\em Proc. 7th Int. Workshop Computational Systems Biology}, pages
  51--54, Luxembourg, Luxembourg, June 2010.

\bibitem{Pfeuty:2011em}
Benjamin Pfeuty, Quentin Thommen, and Marc Lefranc.
\newblock {Robust entrainment of circadian oscillators requires specific phase
  response curves}.
\newblock {\em Biophys. J.}, 100(11):2557--2565, June 2011.

\bibitem{Bunger:2000cs}
Maureen~K Bunger, Lisa~D Wilsbacher, Susan~M Moran, Cynthia Clendenin, Laurel~A
  Radcliffe, John~B Hogenesch, M~Celeste Simon, Joseph~S Takahashi, and
  Christopher~A Bradfield.
\newblock {Mop3 is an essential component of the master circadian pacemaker in
  mammals}.
\newblock {\em Cell}, 103(7):1009--1017, December 2000.

\bibitem{vonGall:2003br}
Charlotte von Gall, Elizabeth Noton, Choogon Lee, and David~R Weaver.
\newblock {Light does not degrade the constitutively expressed BMAL1 protein in
  the mouse suprachiasmatic nucleus}.
\newblock {\em Eur. J. Neurosci.}, 18(1):125--133, July 2003.

\end{thebibliography}

\end{document}